\documentclass[11pt]{amsart}

\usepackage[utf8]{inputenc}
\usepackage{amsmath, amssymb}
\usepackage[all,cmtip]{xy}
\usepackage{booktabs, multirow}
\usepackage{graphicx}
\usepackage{enumerate}
\usepackage{pgf,tikz}
\usetikzlibrary{arrows}

\newcommand{\PP}{\mathbb P}
\newcommand{\CC}{\mathbb C}

\DeclareMathOperator{\Span}{span}

\DeclareMathOperator{\characteristic}{char}
\DeclareMathOperator{\supp}{supp}

\begin{document}

\title{The maximum number of lines lying~on~a~K3~quartic~surface}
\author{Davide Cesare Veniani}

\address{Institut für Algebraische Geometrie \\
Leibniz Universität Hannover \\
Wel\-fen\-garten 1, 30167 Hannover, Germany}
\email{veniani@math.uni-hannover.de}

\date{\today}

\begin{abstract}
We show that there cannot be more than 64 lines on a quartic surface admitting isolated rational double points over an algebraically closed field of characteristic $p \neq 2,\,3$, thus extending Segre--Rams--Schütt theorem. Our proof offers a deeper insight into the triangle-free case and takes advantage of a special configuration of lines, thereby avoiding the technique of the flecnodal divisor. We provide several examples of non-smooth K3 quartic surfaces with many lines.
\end{abstract}

\maketitle
\tableofcontents

\theoremstyle{remark}
\newtheorem{obs}{Observation}[section]
\newtheorem{remark}[obs]{Remark}
\newtheorem{example}[obs]{Example}

\theoremstyle{definition}
\newtheorem{definition}[obs]{Definition} 
\newtheorem{claim}{Claim}

\theoremstyle{plain}
\newtheorem{proposition}[obs]{Proposition}
\newtheorem{theorem}{Theorem}
\newtheorem{lemma}[obs]{Lemma}
\newtheorem{constr}[obs]{Construction}
\newtheorem*{conj}{Conjecture}
\newtheorem{corollary}[obs]{Corollary}

\section*{Introduction}
The aim of this paper is to prove the following result:

\begin{theorem} \label{thm:maintheorem}
 Let $k$ be an algebraically closed field of characteristic $p\geq 0$, $p\neq 2,\,3$ and let $X\subset \PP^3_k$ be a surface of degree 4 over $k$ admitting only isolated rational double points as singularities. Then, $X$ contains at most 64 lines.
\end{theorem}

The problem of finding a bound for the number of lines on quartic surfaces is part of a broader topic, namely the enumerative geometry of lines on surfaces of degree $d$ in projective space. This topic, which has gathered momentum in the last years (see, for example, Boissière--Sarti \cite{boissiere-sarti}, Rams--Schütt \cite{64lines}, Kollár \cite{kollar}), has a long history that dates back to the 19th century. 

The case of smooth surfaces over the complex numbers has drawn so far most attention. Smooth cubic surfaces were thoroughly studied by classical geometers such as Cayley, Clebsch, Salmon, Steiner, Schläfli, Cremona and Sturm. Every smooth cubic surface contains exactly 27 lines which are organized in a highly symmetric way related to the Weyl group of $E_6$.

The general smooth surface of degree $d \geq 4$ contains no lines at all. The first one to state the correct optimal bound of 64 lines for smooth quartic surfaces was B. Segre in 1943 \cite{segre}; nonetheless, his proof contained some major gaps that have been corrected only seventy years later by Rams and Schütt \cite{64lines}. Rams and Schütt used some techniques which were unknown to Segre, above all the theory of elliptic fibrations developed by Kodaira in the 1950's. Segre stated that each line on a smooth quartic surface could meet at most 18 other lines. This was the crucial estimate that Rams and Schütt proved to be false, finding an explicit family of quartics $\mathcal Z$ containing surfaces with a line intersecting 19 or even 20 other lines, which prompted them to work out a new proof. They further examined this family in a follow-up article~\cite{onquartics}.

The problem is still quite open for smooth surfaces of higher degrees. There are some general bounds for $d \geq 5$, but none of them is known to be optimal. Some special cases with particular symmetries have been investigated by Boissière and Sarti~\cite{boissiere-sarti}, where they also find several surfaces with a high number of lines. We refer to their article also for an account of the known bounds.

Although non-smooth cubics had already been classified by Schläfli in the 1860's, it was not until 1979 that the number of lines lying on them was exactly determined by Bruce and Wall~\cite{brucewall}. Non-smooth cubic surfaces always contain less than 27 lines, but one can count the lines with multiplicity -- depending on the number and type of singular points lying on them -- so that the total number is always 27.

Smooth quartic surfaces in $\PP^3$ are K3 surfaces. This paper deals with `mildly' singular quartic surfaces, i.e. with quartic surfaces whose minimal desingularization is a smooth K3 surface; these are precisely the surfaces mentioned in Theorem~\ref{thm:maintheorem}. In this paper we call them `K3 quartic surfaces'; equivalently, K3 quartic surfaces are surfaces in $\PP^3$ of degree~4 admitting only isolated ADE singularities. González Alonso and Rams deal with the case of quartic surfaces with higher singularities in a subsequent paper~\cite{GA-Rams}.

In this paper we extend Rams and Schütt's techniques to our broader setting. The main difficulty still lies in providing a bound for the number of lines meeting a given line. We study the elliptic fibration induced by the given line on the minimal desingularization of the K3 quartic surface; this fibration restricts to a morphism from the strict transform of the line to $\PP^1$. There are two features which make the study of such fibrations much more involved in the K3 case, compared to the smooth one: first, the morphism from the line to $\PP^1$ has always degree 3 if the quartic is smooth, but it has smaller degree as soon as there is a singular point on the line; second, more complicated Kodaira fiber types may appear. Thus, adapting Rams and Schütt's proof to the K3 quartic case forces us to study several new configurations; our results are summarized in Table~\ref{tab:summary}. In addition to the family $\mathcal Z$ found by Rams and Schütt, we discover two new configurations on some non-smooth K3 quartic surfaces in which a line meets more than 18 lines. Explicit examples of such configurations are given in Section~\ref{sec:examples}.

The bound of 64 lines in Theorem~\ref{thm:maintheorem} is sharp and is reached by Schur's quartic, which is smooth. An optimal bound for K3 quartic surfaces with at least one singular point is not known; to our knowledge, the current explicit records are 39 lines over a field of characteristic zero (Example \ref{ex:Gonz-Rams} due to González Alonso and Rams, which is a Delsarte surface) and 48 lines over a field of positive characteristic $\neq 2,\,3$ (Example \ref{ex:48lines}). A K3 quartic surface over $\mathbb C$ with 40 lines exists, but an explicit equation is not known (Example~\ref{ex:implicit}).

A new feature of our proof is that -- unlike Segre and Rams--Schütt -- we do \emph{not} employ a technical tool called `flecnodal divisor' or `flecnodal locus'. Our approach offers a deeper insight into the geometry of the surfaces, and is based on two main ingredients: first, we discover a new, geometrically rich configuration of particular pairs of lines, which we called ``twin lines''; second, thanks to seminal ideas by A. Degtyarev, we take advantage of the well-known lattice theory of K3 surface to tackle the so-called ``triangle-free'' surfaces, i.e. surfaces not containing three lines intersecting pairwise at smooth points.

The newly-developed methods presented here stem from a fruitful synergy between the different points of view of two teams, the one -- S. Rams and M. Schütt -- based in Hannover, Germany, the other -- A. Degtyarev, I. Itenberg and A. S. Sertöz -- mostly in Ankara, Turkey, which for some time worked on the same problem unaware of each other. The Ankara team took a powerful lattice-theoretical approach, which is key to the proof of Lemma~\ref{lemma:(4,6)}. Their work \cite{AIS} has not been published yet.

The proof of our theorem fails over fields of characteristic 2 or 3, mainly because of the presence of quasi-elliptic fibrations; notably, the Fermat quartic surface, considered over an algebraically closed field of characteristic 3, contains exactly 112 lines. Rams and Schütt proved that this is the maximum number that can be achieved in the smooth case \cite{char3}. In the same paper, they also showed a bound of 84 lines for smooth quartics defined over fields of characteristic~2, but the current best example (due to Schütt) contains 68 lines. Throughout the text we highlight the points in which we use that $\characteristic k \neq 2,\,3$.

The paper is organized as follows. In Section~\ref{section:setup} we present some general results on K3 quartic surfaces containing a line. In Section~\ref{sec:valency} we study the number of lines that can intersect a fixed line; the results of this Section are summarized in Table~\ref{tab:summary}. In Section \ref{sec:twinlines} we present the new construction of ``twin lines''. Configuration of lines meeting many other lines are studied in Section~\ref{sec:highvalency}, which also presents the family $\mathcal Z$ found by Rams and Schütt. In Section~\ref{sec:trianglefree} we take care of the triangle-free case. In Section~\ref{sec:prooftheorem} we prove Theorem~\ref{thm:maintheorem}. Finally, in Section~\ref{sec:examples} we list several explicit examples of K3 quartic surfaces with lines meeting many lines, or of K3 quartic surfaces with at least one singular point containing many lines.

\section{Setup} \label{section:setup}

In this Section we collect some basic results and fix the notation which will be used throughout the paper. We always work over an algebraically closed field $k$ of characteristic $p \geq 0$, $p\neq 2,3$. 

\begin{definition} \label{def:K3quartic}
 A \emph{K3 quartic surface} is a surface in $\PP^3$ of degree 4 admitting only isolated rational double points as singularities.
\end{definition}

Any K3 quartic surface $X$ admits a minimal desingularization $\varphi: Z \rightarrow X$, where $Z$ is a smooth K3 surface and $\varphi$ is a birational morphism defined by a numerically effective line bundle 
\[ L:= \varphi^*(\mathcal O_X(1))\] 
of degree~4 (see, for example, \cite[Proposition 1.11]{urabe87}). Given a curve $C$ on $X$, we will always denote by $\hat C$ its strict transform on $Z$.

\begin{lemma}
 Let $X \subset \PP^3$ be K3 quartic surface containing a line $l$, with minimal desingularization $\varphi: Z \rightarrow X$ defined by the nef line bundle $L$. Then, the pencil of planes $\{\Pi_t\}_{t \in \PP^1}$ passing through the line $l$ induces an elliptic fibration on the surface $Z$
 \[
  \pi: Z \rightarrow \PP^1.
 \]
\end{lemma}
\proof The pullbacks by the morphism $\varphi$ of the curves $\Pi_t \cap X$, obtained intersecting $X$ with the pencil of planes containing $l$, define a pencil $\Sigma$ on $Z$. Let $E + \Delta$ be a general member of $\Sigma$, where $\Delta$ is the fixed part; the curve $\hat l$ must be contained in $\Delta$. Given an arbitrary $E' \in |E|$, we have 
$E'+\Delta \in |L|,$
because $|L|$ is complete; therefore, there is a plane $\Pi'$ with $\varphi^* \Pi' = E' + \Delta$. Since 
\[
 \varphi(\supp(E'+\Delta))\supset \varphi(\hat l) = l,
\]
it must be $l \subset \Pi'$, so $E'+\Delta \in \Sigma$; hence $|E|+\Delta = \Sigma$ and, in particular, $\dim |E| = 1$.

The Riemann-Roch formula implies that $E^2 =0$ and $h^1(\mathcal O_Z(E)) = 0$; moreover, $E$ is an irreducible curve of arithmetic genus 1 \cite[Proposition 2.6.]{saint-donat}. The curve $E$ is also smooth: in characteristic zero, this follows from ``generic smoothness'', while for positive characteristic $p \neq 2,\,3$ this can be found, e.g., in Liedtke \cite[Corollary 5.2]{liedtke}.

The base-point-free complete linear system $|E|$ induces a morphism
\[ 
 \pi:= \pi_{|E|}: Z \rightarrow \PP^1;
\]
by what we have just seen, this is an elliptic fibration. \endproof

We denote the restriction of $\pi$ to $\hat l$ again by $\pi$.
\begin{definition}
If the morphism $\pi: \hat l \rightarrow \PP^1$ is constant, we say that $l$ has degree 0; otherwise, the \emph{degree} of $l$ is the degree of the morphism $\pi: \hat l \rightarrow \PP^1$.
\end{definition}

\begin{definition}
The \emph{singularity} of a line $l$ is the number of singular points lying on $l$.
\end{definition}

\begin{proposition} \label{prop:degree}
The degree $d$ of a line $l \subset X$ is given by 3 minus the intersection multiplicities at each singular point on the line $l$ of a general residual cubic in the pencil and the line $l$:
  \begin{equation} \label{eq:degree}
   d = 3 - \sum_{P \in l \cap S} I_P(E_t, l),
  \end{equation}
  where $S\subset X$ is the set of singular points of $X$, $E_t$ is the residual cubic in the plane $\Pi_t$, for a general $t \in \PP^1$, and $I_P(E_t, l)$ is the intersection multiplicity of $E_t$ and $l$ at $P$.
  In particular,
  \[
      d \leq 3 - s, \quad \text{and} \quad d = 3 \iff s = 0,
  \]
  where $s$ is the singularity of $l$.
  The morphism $\pi:\hat l \rightarrow \PP^1$ is a separable map outside characteristic 2 and 3.
\end{proposition}
\proof We will presently see that the integer
    \[
    d:= \deg(\pi) = E.\hat l 
    \]
can vary between $0$ and $3$, by computing it explicitly; the morphism $\pi:\hat l \rightarrow \PP^1$ is thus a separable map outside characteristic 2 and 3.

Suppose that the coordinates in $\PP^3$ are $x_0,x_1,x_2,x_3$ and that the quartic $X$ is given by
\begin{equation} \label{eq:surfaceX}
 X: \sum a_{i_0 i_1 i_2 i_3} x_{0}^{i_0} x_{1}^{i_1} x_{2}^{i_2} x_{3}^{i_3} = 0, \quad a_{i_0 i_1 0 0} = 0 \text{ for all $i_0$, $i_1$},
\end{equation}
where the sum is taken over all quadruples of nonnegative integers which sum to 4, and that we have chosen the coordinates so that the line $l$ on $X$ is given by $x_0 = x_1 = 0$ (whence the condition on the $a_{i_0 i_1 0 0}$).

Let us rewrite the equation of $X$ like this:
\[
 X: x_0 \alpha(x_2, x_3) + x_1 \beta(x_2,x_3) + \text{terms containing $x_0^2$, $x_0 x_1$ or $x_1^2$}.
\]
The forms $\alpha$ and $\beta$ have degree 3 and we write them down explicitly:
\begin{align} \label{eq:alphaebeta}
\begin{split}
 \alpha(x_2, x_3) & = a_{1030} x_2^3 + a_{1021} x_2^2 x_3 + a_{1012} x_2 x_3^2 + a_{1003} x_3^3, \\
 \beta(x_2, x_3) & = a_{0130} x_2^3 + a_{0121} x_2^2 x_3 + a_{0112} x_2 x_3^2 + a_{0103} x_3^3.
\end{split}
\end{align}

The intersection of the quartic with the pencil of planes $x_0 = t x_1$ is given by the line $l$ and the residual cubics $E_t$. One sees explicitly that the intersection of $E_t$ with $l$ is given by the points $[0:0:x_2:x_3]$ satisfying
\begin{equation} \label{eq:penciloncurve}
 t \alpha(x_2, x_3) + \beta(x_2, x_3) = 0.
\end{equation}

In particular, the degree of $\pi:\hat l \rightarrow \PP^1$ is given by 3 minus the common roots of $\alpha$ and $\beta$ counted with multiplicity (note that $\alpha$ and $\beta$ cannot be identically zero at the same time, otherwise $l$ would be a line of singular points). Observe that $l$ contains a singularity at the point
\[ P = [0:0:x_2:x_3] \] 
exactly when $[x_2:x_3]$ is a common root of $\alpha$ and $\beta$; moreover, the multiplicity of this common root is exactly $I_P(E_t,l)$ for a general $t \in \PP^1$. This proves formula \eqref{eq:degree}; in particular, if $\alpha$ and $\beta$ have no roots in common or, equivalently, if there are no singularities on $l$ -- and only in that case -- the degree of the associated morphism $\pi: \hat l \rightarrow \PP^1$ is 3.

If $\alpha$ and $\beta$ have all roots in common, i.e. they are multiple of each other, then there is exactly one plane whose intersection with the quartic contains $l$ as a non-reduced component (if there were more, the surface would have worse singularities than isolated rational double points): this plane is the only plane tangent to the surface along the line $l$. This means precisely that the degree of the morphism $\pi: \hat l \rightarrow \PP^1$ is zero, that is to say, $\hat l$ is a component of a fiber of the morphism $\pi: Z \rightarrow \PP^1$.
\endproof

Let now $X$ be any K3 quartic surface (not necessarily containing a line). If $X$ has a singularity in $P$ we can choose coordinates so that $P = [0:0:0:1]$ and the defining equation of $X$ is
\begin{equation} \label{eq:f2f3f4} 
 x_3^2 f_2(x_0, x_1, x_2) + x_3 f_3(x_0, x_1, x_2) + f_4(x_0, x_1, x_2) = 0,
\end{equation}
where the $f_i$'s are homogeneous forms of degree $i$.

\begin{definition}
We call these forms the \emph{(second, third, fourth) Taylor coefficients} of $X$ at $P$. 
\end{definition}

Since we are considering only rational double points, the form $f_2$ is not identically zero and the equation $f_2 = 0$ defines the \emph{tangent cone} of $X$ at $P$. 

\begin{lemma} \label{lemma:linesthroughsingularpoint}
 Suppose $P$ is a singular point of a K3 quartic surface $X$. Then there are at most 8 lines contained in the surface $X$ passing through it.

 Moreover, if there are more than 6, then the second and the third Taylor coefficients of $X$ at $P$ share a common factor. If there are 8, then the second Taylor coefficient of $X$ at $P$ must divide the third.
\end{lemma}
\proof Consider equation \eqref{eq:f2f3f4}. A line parametrized by $t \mapsto [a t: b t: c t:1]$ is contained in $X$ if and only if 
$[a:b:c]$ is a point of intersection of the three plane curves of degree $i=2,3,4$ defined by
\[ f_i = 0. \]
Recall that by B\'ezout's theorem two plane curves of degree $d$ and $e$ without irreducible components in common have at most $d\cdot e$ distinct point in common.

Note first that $f_2$, $f_3$ and $f_4$ cannot all have a common irreducible component, otherwise the surface $X$ would be reducible.

Suppose first that $f_2$ is irreducible. Then, by B\'ezout's theorem, it intersects $f_3$ in at most 6 points, unless it is an irreducible component $f_3$; in this case $f_2$ divides $f_3$ and, since $f_2$ is not a component of $f_4$, $f_2$ and $f_4$ have at most 8 common points, which is what we claimed.

Suppose now that $f_2 = gh$ is the union of two lines $g= 0$, $h=0$, which may be identical or different. If none of these lines is a component of $f_3$, then the number of common solutions of $f_2$ and $f_3$ is at most 6. 

Hence, if the number of common solutions is bigger than 6, the curves defined by $f_2$ and $f_3$ have at least one common irreducible component, say $g$. Since each common component of $f_2$ and $f_3$ is not a component of $f_4$, $g$ gives at most 4 solutions with $f_4$. If $h$ is not a component of $f_3$, then they intersect in at most 3 distinct point, so the number of intersection points is at most $3+4=7$. Therefore, in order to have 8 distinct solutions the two lines $g,h$ must be different and also $h$ must be a component of $f_3$, which implies that the polynomial $f_2$ divides $f_3$. \endproof

\section{Valency of a line} \label{sec:valency}

In this section, $X$ will always be a K3 quartic surface in $\mathbb P^3$ and $\varphi:Z\rightarrow X$ will denote its minimal desingularization defined by a numerically effective line bundle $L$ of degree~4. We are interested in finding a bound on the number of lines meeting a given line $l\subset X$. We need to distinguish between the lines that meet $l$ at smooth points and those that meet it at singular points.

\begin{definition}
    A \emph{line} on $Z$ is a smooth rational curve $C \subset Z$ such that $C\cdot L = 1$.
\end{definition}

A smooth rational curve $C \subset Z$ is a line if and only if it is the strict transform of a line on $X$.

\begin{definition}
Given a line $l \subset X$, the \emph{valency} of $l$, denoted by $v(l)$, is the number of lines $m$ on $Z$ that intersect the strict transform $\hat l$ of $l$; in other words, the valency of $l$ is the number of lines $m$ on $X$ such that their strict transform $\hat m \subset Z$ intersects the strict transform $\hat l$ of $l$.
\end{definition}

\begin{definition}
Given a line $l$ on K3 quartic surface $X$, the \emph{extended valency} of $l$, denoted by $\tilde v(l)$, is the number of lines on $X$ that intersect $l$.
\end{definition}

\begin{definition}
    A plane containing $l$ such that the corresponding residual cubic splits into three lines (not necessarily distinct) is called a $p$-fiber of $l$; a plane containing $l$ such that the corresponding residual cubic splits into a line and an irreducible conic is called a $q$-fiber of $l$.
\end{definition}

\begin{definition}
    We say that a line $l$ is of type $(p,q)$, for $p,q\geq 0$, if the number of planes containing $l$ and splitting into three lines (resp. a line and an irreducible conic) is $p$ (resp. $q$).
\end{definition}

For any line $l \subset X$ of positive degree and type $(p,q)$, the following relations are obvious:
\begin{equation} \label{eq:vl}
    v(l) \leq \tilde v(l) = 3\,p + q.
\end{equation}
Equality holds if $l$ has degree 3. On the other hand, lines of degree 0 behave in a very special way, as the following lemma shows.

\begin{lemma}
If $l \subset X$ is a line of degree $0$, then $v(l) \leq 2$.
\end{lemma}
\proof A line $l\subset X$ is of degree 0 if and only if there exists a plane $\Pi$ which is tangent to $X$ along $l$, i.e. $l$ appears with multiplicity at least 2 in the intersection of $\Pi$ with $X$. The residual conic on this plane might split into two lines. On the other hand, all other lines meeting $l$ must pass through one of the singular points of $l$, thus not contributing to the valency of $l$. \endproof

Given any line $l \subset X$, we can compute the Euler number of $X$ using the elliptic fibration induced by $l$ as the sum of the Euler number of the singular fibers. Observing that the contribution of any $p$-fibers to this sum is at least 3, and the contribution of any $q$-fiber is at least 2, we get
\begin{equation} \label{eq:3p+2q}
    3\,p + 2\,q \leq 24.
\end{equation}
In particular, it is immediate to note the following lemma.
\begin{lemma} \label{lemma:no-p-fibers}
    If $l \subset X$ has no $p$-fibers, then $\tilde v(l) \leq 12$.
\end{lemma}

An immediate consequence of the inequalities \eqref{eq:vl} and \eqref{eq:3p+2q} is that $v(l)$ and $\tilde v(l)$ cannot be greater than 24. This bound is not sharp. In this section we will prove that $v(l) \leq 20$ for any line $l \subset X$ (see Table~\ref{tab:summary}); it is also possible to prove that $\tilde v(l)\leq 20$ (see Section~\ref{sec:examples}).

A first useful technique to find bounds for $v(l)$ for a line $l\subset X$ of positive degree is to find the points of intersection of the residual cubics $E_t$ and $l$ which are inflection points for $E_t$. By ``inflection point'' we mean here a point which is also a zero of the hessian of the cubic. In fact, if a residual cubic $E_t$ contains a line as a component, all the points of the line will be inflection points of $E_t$. Moreover, if $E_t$ contains an irreducible conic $C$ and a line $m$ s components, then the only points of $C$ which are inflection points are those that also lie on $m$.

Supposing that the surface $X$ is given as in equation \eqref{eq:surfaceX} and $l$ by $x_0 = x_1 = 0$, we write out the hessian of the equation defining $E_t$ and then we substitute $x_1 = 0$, obtaining the condition
\begin{equation} \label{eq:hessian}
 \left| \begin{array}{ccc}
\zeta_3  & \eta_2   & \theta_2 \\
\eta_2   & \iota_1  & \kappa_1 \\
\theta_2 & \kappa_1  & \lambda_1 \end{array} \right| = 0.
\end{equation}
Here $\zeta_3,\, \eta_2,\, \theta_2,\, \iota_1,\, \kappa_1,\, \lambda_1$ are polynomials in $t$ with degree equal to the respective indices and with forms of degree 1 in $(x_2,x_3)$ as coefficients. So this condition is actually given by a polynomial of degree 5 in $t$, the coefficients of which are forms of degree 3 in $(x_2,x_3)$. 

We want now to find the number of lines intersecting $l$ by studying the common solutions of \eqref{eq:penciloncurve} and \eqref{eq:hessian}. It is convenient to extend Segre's nomenclature~\cite{segre}.

\begin{definition}
 A line $l$ of degree $d \geq 1$ is a \emph{line of the first kind} if when substituting $t$ from \eqref{eq:penciloncurve} into \eqref{eq:hessian} we do \emph{not} get the zero form. Otherwise, we will say that $l$ is a \emph{line of the second kind}.
\end{definition}

In the case of a line of the first kind, the zeros of the form obtained by eliminating $t$ correspond to points of the line $l$ that are inflection points for some residual cubic. Note that if 2 (respectively 3) lines intersect $l$ at the same point, then the corresponding zero will have multiplicity greater or equal to 2 (respectively 3): this can be checked by a local computation. This observation yields the following lemma.

\begin{lemma} \label{lemma:1stkind}
Let $l \subset X$ be a line of degree $d \geq 1$ of the first kind; then,
\[
    v(l) \leq 3 + 5\,d.
\]
\end{lemma}

Finding a bound for a line of the second kind is more involved. We divide our analysis according to the degree of $l$. 

We first need a technical lemma. We recall that the Mordell-Weil group of an elliptic surface, i.e. the group of sections of the elliptic fibration, acts on each smooth fiber by translation, and that this action extends to the singular fibers \cite{miranda}.

\begin{lemma} \label{lemma:eulernumber}
    Let $Y \rightarrow \Gamma$ be an elliptic surface over an algebraically closed field $k$ and suppose that the fibration is endowed with an $n$-torsion section $\sigma$, with $\characteristic(k) \nmid n$. Then, the minimal desingularization of the quotient $Y/G$ by the group $G$ generated by the action of $\sigma$ is an elliptic surface $Y' \rightarrow \Gamma$ such that
\[
 e(Y) = e(Y').
\]
\end{lemma}
\proof[Sketch of proof]
     Let us denote by $Y'$ the minimal desingularization of $Y/G$. Since $\sigma$ acts fiberwise, $S'$ is also an elliptic surface over $\Gamma$.  Let $f: Y \dashrightarrow Y'$ be the rational map induced by composition; it corresponds to a morphism $f: Y \setminus F \rightarrow Y'$, where F is a finite set.

    If $\omega$ is a regular 1-form, its pullback $f^*\omega$ is a rational 1-form and is regular on $S \setminus F$. The pullback $f^*\omega$ is not the zero form, since $\characteristic(k) \nmid n$. Since the poles of a non-zero differential form are divisors, $f^*\omega$ is regular on all $Y$, hence there is an injective map $f^*: \Gamma(Y',\Omega_{Y'}) \rightarrow \Gamma(Y,\Omega_Y)$. In particular we have $q(Y') \leq q(Y)$. The same applies to 2-forms, so $p_g(Y') \leq p_g(Y)$.

    Then $f$ induces an isogeny on the generic fibers and the dual isogeny induces a rational map $Y' \dashrightarrow Y$. The same argument works, so we have the equalities $q(Y) = q(S')$ and $p_g(Y) = p_g(Y')$. Since the surfaces are elliptic, they both have $K^2 = 0$. We conclude by applying Noether's formula.
\endproof

\begin{proposition} \label{prop:degree3}
If $l \subset X$ is a line of the second kind of degree $3$, then $v(l) \leq 20$. Moreover, if $v(l) >16 $, then the morphism $\pi: \hat l \rightarrow \PP^1$ has exactly two points of ramification. If $v(l) = 19$, then the line has type $(p,q) = (6,1)$ and the $q$-fiber is a ramified fiber of type $I_n$, $n\geq 2$; if $v(l) = 20$, then the line has type $(p,q) = (6,2)$ and both $q$-fibers are ramified fibers of type $I_n$, $n\geq 2$.
\end{proposition}
\proof Since $\pi: \hat l \rightarrow \PP^1$ has degree 3, by the Riemann-Hurwitz formula it can have 4, 3 or 2 points of ramification, of ramification index respectively $(2,2,2,2)$, $(3,2,2)$ or $(3,3)$. For the sake of simplicity, we call these cases A, B and C.

The morphism $\pi$ corresponds to a field extension $k(\hat l)\supset k(\PP^1)$ of degree~3. This extension is Galois if and only if we are in case C. In fact, the index of ramification at a point $P \in \hat l$ is equal to the order of the inertia group of the corresponding place in $k(\hat l)$. Since the inertia group is a subgroup of the Galois group, its size must divide the size of the Galois group; hence, this extension can never be a Galois extension if there is a point of ramification index 2.
In all three cases, we can draw the following commuting diagram:
\[
 \xymatrix{
 W \ar[d] \ar@{-->}[r] & Z \ar[d]^\pi \\
 \hat l \ar[r]^{\pi} & \PP^1}
\]
where $W$ is the minimal desingularization of the surface $Z \times_{\PP^1} \hat l$ and the map $W \dashrightarrow Z$ is a rational dominant map of degree 3.

In the non-Galois cases A and B, the Galois closure of the extension $k(\hat l)\supset k(\PP^1)$ corresponds to a morphism $\psi: \Gamma \rightarrow \hat l$ of degree 2, where $\Gamma$ is a smooth curve of genus 1 or 0, respectively. We need to perform a further base change, obtaining the following commuting diagram:
\[
 \xymatrix{
Y \ar[d] \ar@{-->}[r] & W \ar[d] \ar@{-->}[r] & Z \ar[d]^\pi \\
\Gamma \ar[r]^{\psi} & \hat l \ar[r]^{\pi} & \PP^1}
\]
where $Y$ is the minimal desingularization of the surface $W \times_{\hat l} \Gamma$ and dashed arrows represent rational dominant maps. In the Galois case C, we can perform a trivial base change, obtaining the same diagram with $\Gamma = \hat l$ and $Y = W$. We set $\eta := \pi \circ \psi$ and $d = \deg \eta$; note that $d = 6$ in cases A and B and $d = 3$ in case C.

The inclusion $\hat l \hookrightarrow Z$ lifts to a section $\hat l \rightarrow W$, which in turn lifts to a section $s_0: \Gamma \rightarrow Y$. The Galois action provides us with two more sections $s_1, s_2: \Gamma \rightarrow Y$. We will call these three sections `Galois sections'. Since the line $l$ is of the second kind and since the three inflection points that $l$ cuts on the general cubic are of course allineated, we can choose $s_0$ to be the 0-section in the Mordell-Weil group of the elliptic surface $Y \rightarrow \Gamma$, so that $s_1$ and $s_2$ become 3-torsion sections inverse to each other.

The action of the Galois sections induces a rational map $Y \dashrightarrow Y'$, where $Y$ and $Y'$ are two elliptic surfaces over $\Gamma$ with the same Euler number, by Lemma \ref{lemma:eulernumber}. To each singular fiber $G$ of $Y$ there corresponds a singular fiber $G'$ of $Y'$ and the type of $G$ is determined univocally by the type of $G$ and by which components of $G$ are met by the torsion sections.

Let $F = \pi^{-1}(t)$ be a singular fiber of $\pi: Z \rightarrow \PP^1$ and let $\Pi \supset l$ be the corresponding plane in $\PP^3$. Let $E$ be the residual cubic in the plane $\Pi$. With a slight abuse of language, we will say that $F$ is `unramified', `ramified of index 2' or `ramified of index 3', according to whether the configuration of $\eta^{-1}(t)$ consists of
\begin{enumerate}
 \item $d$ distinct points ($t$ is not a branch point);
 \item 3 points of ramification index 2;
 \item $d/3$ points of ramification index 3.
\end{enumerate}
Observe that the ramified fibers are in one-to-one correspondence with the ramification points of $\pi: \hat l \rightarrow \PP^1$ and have the same ramification indices.

\emph{§ Unramified fibers.} Suppose $\eta^{-1}(t)$ consists of $n$ distinct points. In this case the fiber $F$ is replaced by $n$ fibers on $Y$. Choose one of them and call it $G$. Note that $F$ and $G$ are of the same type. Since $G$ accommodates 3-torsion sections, $G$ (hence also $F$) must be of type $I_n$ $(n\geq 0)$, $IV$ or $IV^*$ in Kodaira's notation \cite[Table (VII.3.4)]{miranda}.

Suppose that $F$ and $G$ are of type $I_n$ and suppose that the sections $s_0$, $s_1$ and $s_2$ meet the same component of $G$. Recalling that the sections $s_i$ are induced by the strict transform of $l$ in $Z$, the former case happens if and only if the line $l$ meets only one component of $E$, i.e. $E$ is irreducible and gives no contribution to the number of lines meeting $l$). $G$ must correspond to a fiber $G'$ of type $I_{3n}$ on $Z_2'$.

If $F$ and $G$ are of type $I_n$, but the sections $s_i$ intersect different irreducible components, then $n$ must be a multiple of 3 and $G$ corresponds to a fiber $G'$ of type $I_n$ on $Z_2'$. The residual cubic $E$ splits into three lines.

If $F$ and $G$ are of type $IV$ and $IV^*$ then the sections $s_i$ must meet different components of $G$, hence the residual cubic must fully split and the fiber $G'$ on $Z_2'$ corresponding to $G$ has the same type of $F$ and $G$.

\emph{§ Ramified fibers of index 2.} Suppose $\eta^{-1}(t)$ consists of 3 points of ramification index 2. This is only possible in cases A and B, i.e. when $d = 6$. In this case $F$ is replaced by three fibers on $Y$, whose type can be read off from \cite[p. 64, Table (VI.4.1)]{miranda}. A priori, the fiber $F$ can be of type $I_n$, $I_n^*$, $II$, $IV$, $II^*$ or $IV^*$, yielding three fibers on $Y$ of type $I_{2n}$, $I_{2n}$, $IV$, $IV^*$, $IV^*$ or $IV$ respectively, all of which could accommodate 3-torsion sections.

    We can exclude fibers of type $I_n$, though. Indeed, since the the line $l$ meets the residual cubic in $\Pi$ at inflection points, $l$ cannot be tangent to the residual cubic, otherwise it would have intersection of order 3 and $\Pi$ would not correspond to a ramification point of $\pi$ of index 2. Hence, $\hat l$ cannot be tangent to the fiber $F$ (since all blowups of the desingularization happen outside of $l$) and, therefore, $\hat l$ meets $F$ in a node. However, on each new fiber on $Z_2$, two of the Galois sections, say $s_0$ and $s_1$, meet the same component, and the third one $s_2$ meets a different component: this is impossible since we could choose $s_0$ to be the 0-section, but $s_1$ and $s_2$ could not be the inverse of each other.

    We also deduce that the residual cubic $E$ in the plane $\Pi$ corresponding to $F$ cannot split into three different lines. Indeed, these lines could not be concurrent since $l$ would pass through the node and ramification of index 3 would occur. But if they were not concurrent, then they would form a `triangle'; after blowing up the singular points, the fiber $F$ would still contain a `cycle', hence it should be of type $I_n$ ($n\geq 3$), which we have just ruled out.

    The residual cubic $E$ cannot split into a line and a conic, because if the two where secant, then they would form a `cycle' (of length 2) and this cycle would lead to a fiber of type $I_n$, while if they where tangent then $l$ would pass through the point of tangency (since it cannot be tangent to the conic) and we would get a fiber of type $III$, which is also excluded; moreover, the residual cubic $E$ cannot be an irreducible cubic with a node, since $l$ should pass through the node and we would have a fiber of type $I_1$ (recall that all points on $l$ are smooth because $\pi:Z\rightarrow \PP^1$ has degree 3); finally, $E$ cannot be a triple line, otherwise ramification of index 3 would occur.

    Hence, we are left with very few possibilities: either $E$ is an irreducible cubic with a cusp, $l$ passes through the cusp and we have a fiber of type $II$, or $E$ splits into a double line and another line, hence the fiber $F$ contains a component of multiplicity 2. This rules out a fiber $F$ of type $IV$, too.

\emph{§ Ramified fibers of index 3.}
Suppose $\eta^{-1}(t)$ consists of $d/3$ points of ramification index 3. In this case, $F$ is replaced by $d/3$ fibers on $Y$. As before, we be read off their fiber type from \cite[p. 64, Table (VI.4.1)]{miranda}: the fiber $F$ can be of type $I_n$, $IV$ or $IV^*$.

If the residual cubic in $\Pi$ has a non-reduced component, then it must lead to a fiber of type $IV^*$. If the residual cubic is composed of three distinct lines, then, in order to have ramification of type 3, they must be concurrent and the line $l$ must pass through the node; thus, there cannot be singular points of the surface on the three lines (since this would result in a fiber outside Kodaira's classification) and the fiber must be of type $IV$.

If the fiber $F$ is of type $I_1$, then the residual cubic must be irreducible; hence, it gives no contribution to the lines meeting $l$. If the fiber $F$ is of type $I_n$, $n\geq 2$, then the residual cubic splits into a line plus a conic (it cannot split into three lines, otherwise we could not have ramification of index 3). In each case, the three Galois sections must meet the same component on each of the two fibers of type $I_{3n}$ on $Z_2$ (this component comes from the node of the residual cubic through which $l$ passes); therefore, we get two fibers of type $I_{9n}$ on $Z_2'$.

Three concurrent lines correspond to a fiber $F$ of type $IV$. A double or a triple line  must lead to a fiber $F$ of type $IV^*$.

\emph{§ Conclusions.} Table \ref{tab:degree3} summarizes what we have proven so far, where the column ``difference'' represents the contribution to the difference of the Euler numbers of $Z_2$ and $Z_2'$ due to the fibers obtained from $F$ by base change and the column ``$v_t(l)$'' stands for the number of lines contained in the corresponding plane $\Pi_t$; of course,
\[
v(l) = \sum_{t\in \PP^1} v_t(l).
\]

\begin{table}[t]
\caption{\small Possible singular fibers for a line of the second kind of degree 3.}
\label{tab:degree3}
\begin{tabular}{llllll} \toprule
ramification & fiber on $Z$ & fibers on $Y$      & fibers on $Y'$     & difference & $v_t(l)$\\
\midrule
\multirow{4}{*}{\small unramified} & $I_{n}$      & $d \times I_{n} $ & $d \times I_{3n}$ & $+12\,n$   & 0 \\
                                   & $I_{3n}$     & $d \times I_{3n}$ & $d     \times I_n$    & $-12\,n$   & 3 \\
                                   & $IV$         & $d \times IV$     & $d \times IV$     & 0          & 3 \\
& $IV^*$       & $d \times IV^*$   & $d \times IV^*$   & 0          & 3 \\
  \midrule
\multirow{5}{*}{\parbox[t]{1.6cm}{\small ramification\\of index 2}} & $I^*_{n}$ & $3 \times I_{2n} $ & $3 \times I_{6n}$ & $+12\,n$ & 2 \\
  & $I^*_{3n}$ & $3 \times I_{6n} $ & $3 \times I_{2n}$ & $-12\,n$ & 2 \\
  & $II$ & $3\times IV$ & $3\times IV$ & 0 & 0\\
  & $II^*$ & $3 \times IV^*$ & $3 \times IV^*$ & 0 & 2\\
  & $IV^*$ & $3 \times IV$ & $3 \times IV$ & 0 & 2\\
 \midrule
\multirow{4}{*}{\parbox[t]{1.6cm}{\small ramification\\of index 3 \\ ($d' = d/3$)}} & $I_1$ & $d' \times I_{3}$ & $d' \times I_{9}$ & $+6\,d'$ & 0 \\
 & $I_n$ ($n\geq 2$)& $d' \times I_{3n}$ & $d' \times I_{9n}$ & $+6\,nd'$ & 1\\
 & $IV$ & $d' \times I_0$ & $d' \times I_0$ & 0 & 3\\
 & $IV^*$ & $d'\times I_0$ & $d' \times I_0$ & 0 & $\leq 2$\\
 \bottomrule
\end{tabular}
\end{table}

 By Lemma~\ref{lemma:eulernumber} the Euler numbers of $Y$ and $Y'$ must balance out. In case A, there are always four ramified fibers, so their contribution to the Euler number is always at least 8. Considering the possible combinations of fibers, one can see that each time we get $3 \,n$ lines we must pay with a contribution of at least $4\, n$ to the Euler number, so the number of lines intersecting $l$ is not greater than 12.

In case B, the contribution to the Euler number coming from the ramified fibers of index 2 is at least 4, without any contribution to the number of lines. Again, looking at the possible combinations, one can see that we need a further contribution of at least $4\,n$ to the Euler number each time we get $3\,n$ lines, except when we have a ramified fiber of type $I_{n}$, $n\geq 2$, (there can be at most one) paired with $n$ unramified fibers of type $I_3$, in which case we get $3\,n +1$ lines for a loss of $4n$ in the Euler number. Hence, the maximal number of lines meeting $l$ is 16.

In case C, A direct inspection of the possible combinations yields a bound of 20 lines meeting $l$. The line $l$ can meet 19 or 20 lines only if there are one or two ramified fibers of type $I_n$, $n\geq 2$.
\endproof

\begin{proposition} \label{prop:degree2}
If $l \subset X$ is a line of the second kind of degree $2$, then $v(l) \leq 10$.
\end{proposition}
\proof Referring to the notation in Section \ref{section:setup}, in this case there must be a singular point $P$ on the line $l$ and it must correspond to a simple common root of $\alpha$ and $\beta$, so that the generical residual cubic intersects $l$ at $P$ with multiplicity 1.

Having degree 2, the morphism $\pi: \hat l \rightarrow \PP^1$, ramifies in two distinct points of ramification index 2, by the Riemann-Hurwitz formula. By a slight abuse of language we will say, for example, that ``the point $Q \in l$ is a ramification point'' instead of saying that ``the only point $\hat Q \in \hat l$ that maps to $Q$ through $\varphi$ is a ramification point''.

Up to projective equivalence, we can suppose that $P = [0:0:0:1]$ is the singular point on $l$, which is the same as requiring $x_2$ to be the common root of $\alpha$ and $\beta$; this means that in equation \eqref{eq:surfaceX} we have
\begin{equation} \label{eq:caseB1}
  a_{0103} = a_{1003} = 0.
\end{equation}

In addition, after a suitable change of coordinates we can suppose that
\begin{equation} \label{eq:caseB2}
 a_{0112} = a_{1030} = a_{1021} = 0 \quad \text{and} \quad a_{1012} = a_{0130} = 1.
\end{equation}
In fact, we can assume that $Q = [0:0:1:0]$ is a point of ramification relative to the plane $\Pi_1: x_1 = 0$ and that the residual cubic in $\Pi_0: x_0 = 0$ has at least double intersection with $l$ in $P$. Note that $P$ is a point of ramification if and only if $a_{0121} = 0$. The coefficient $a_{1012}$ and $a_{0130}$ must be different from 0, so we can normalize them to 1.

By eliminating $t$ from \eqref{eq:penciloncurve} and \eqref{eq:hessian}, we get a homogeneous form of degree $5\cdot 2 +3 = 13$ in $w,\,z$, which we denote by $h_{13}$. Since the line $l$ is of the second kind, this form vanishes. By looking at the coefficients of $x_3^{13}$ and $x_2 x_3^{12}$ one finds that $a_{0202}$ must vanish; hence, the tangent cone at $P$, which is given by
\[
 f_2 : {\left(a_{2002} x_{0} + a_{1102} x_{1} + x_{2}\right)} x_{0} = 0,
\]
is the union of two distinct planes, whose intersection is a line different from $l$. This means that the point $P$ can be neither of type $A_1$ (since the tangent cone would be irreducible) nor of type $D_n$ nor of type $E_n$ (since the tangent cone would be a double plane). Therefore, $P$ is of type $A_n$, with $n\geq 2$, and on the minimal desingularization $Z$ we have $n$ exceptional smooth rational exceptional divisor $\Delta_1, \ldots, \Delta_n$, such that $\Delta_i.\Delta_{i+1} = 1$ for $i=1,\ldots,n-1$, and $\Delta_i.\Delta_j = 0$ otherwise (as long as $i\neq j$).

The fact that the line of intersection of the two planes of the tangent cone is different from $l$ tells us that the strict transform $\hat l$ of $l$ meets one `extremal' exceptional component, say $\Delta_1$, while the strict transform of a general residual cubic meets the other `extremal' exceptional component $\Delta_n$; in fact, the two `extremal' components parametrize the tangent directions in the two planes of the tangent cone.

Consider now the following commuting diagram
\begin{equation} \label{eq:basechange}
\xymatrix{
 W \ar[d] \ar[r] & Z \ar[d]^\pi \\
 \hat l \ar[r]^{\pi} & \PP^1}
\end{equation}
where $W$ is the minimal desingularization of $Z \times_{\PP^1} \hat l$.

Note that the field extension $k(\hat l) \supset k(\PP^1)$ corresponding to $\pi$ has degree 2; hence, it is always Galois outside characteristic 2. Therefore, we have three sections $s_0,\,s_1,\,s_2:\hat l \rightarrow W$, which we will call `Galois sections'. We can choose one of them to be the zero section; since $l$ is of the second kind, the other two are 3-torsion sections. Observe that two sections map one-to-one onto $\hat l$ through $\psi$, whereas the third section maps two-to-one onto $\Delta_n$. The Galois sections induce a rational map $W \dashrightarrow W'$, where $e(W) = e(W')$, as explained in Lemma~\ref{lemma:eulernumber}.

Let $F_t = \pi^{-1}(t)$ be a singular fiber of the elliptic fibration $\pi: Z \rightarrow \PP^1$ induced by $l$ corresponding to a plane $\Pi$.

\emph{§ Unramified fibers.} If $t$ is not a branch point of  $\pi: \hat l\rightarrow \PP^1$, then $F_t$ has type $I_n$ $(n\geq 1)$, $IV$ or $IV^*$, since on $W$ it is substituted by two fibers of the same type and these must accommodate 3-torsion.

Suppose $F_t$ is a fiber of type $I_n$. If the residual cubic in $\Pi$ is irreducible, then the three Galois sections meet the same component; hence, we get two fibers of type $I_{3n}$ on $W'$ and these fibers do not contribute to the valency of $l$. On the other hand, if the residual cubic in $\Pi$ is reducible, then $n$ must be divisible by 3 and we get two fibers of type $I_{3m}$ on $W$ and two of type $I_m$ on $W'$, where $n = 3\,m$.

\emph{§ Ramified fibers.} If $t$ is a branch point of $\pi: \hat l\rightarrow \PP^1$, then a priori $F$ can have type $I_n$ $(n\geq 1)$, $I_n^*$ $(n\geq 1)$, $II$, $IV$, $II^*$ or $IV^*$ (again, see \cite[p. 64, Table (VI.4.1)]{miranda}). We can exclude type $I_n$ and $IV$, though.

We call $\hat P$ the point of intersection of $\Delta_1$ with $\hat l$. There exists exactly one fiber $F_0$ containing $\hat P$; let us denote by $\Pi_0$ the corresponding plane (in our parametrization $\Pi_0$ is given by $x_0 = 0$). Note that the fiber $F_0$ must contain $\Delta_1,\ldots,\Delta_{n-1}$ as irreducible components plus the strict transform of the components of the residual cubic $E_0$ in $\Pi_0$.

If $F_0$ is a ramified fiber, one can see by a local computation that the residual cubic $E_0$ must split into three lines passing through $P$: in fact, setting $a_{0121} = 0$ (which was the condition for ramification in $P$) the residual cubic in $x_{0} = 0$ has no term containing $x_3$. The three lines can be all distinct or they might coincide. In any case, we have no cycles and more than three components; hence, we can exclude type both type $I_n$ aand $IV$.

Suppose now that $F$ is a ramified fiber different from $F_0$. The corresponding residual cubic $E$ has thus intersection multiplicity 1 with $l$ at $P$ and 2 at another point $Q \in l$. $P$ is the only singular point of $X$ on $l$ by Proposition~\ref{prop:degree}, so $Q$ must be a smooth point of $X$. Moreover, since $Q$ is an inflection point of $E$ because $l$ is of the second type, $E$ and $l$ cannot meet tangentially in $Q$, otherwise the intersection multiplicity would be 3.

Hence, if the cubic $E$ is irreducible, then $Q$ must be a cusp, and $F$ is of type $II$. In fact, if $Q$ were a node, then two of the Galois sections on $Z_1$ would meet the same component of the resulting $I_2$-fiber on $Z_1$ and the third would meet a different one, which is impossible.

The cubic $E$ cannot split into a line and a conic, because in this case $Q$ would be a point of intersection of the line and the conic, giving rise either to a fiber of type $III$ (which we excluded a priori) or to a fiber of type $I_n$ ($I_2$ if there are no surface singularities in the plane relative to $E$, otherwise $I_n$ with $n>2$) with an impossible configurations of torsion sections as before.

If the cubic $E$ split into three distinct lines, they could not be concurrent because $F$ is a ramified fiber different from $F_0$, so again this would lead to an impossible configuration of torsion sections. Finally, if $E$ splits into three lines not all distinct, then $F$ contains a nonreduced component, so fiber types $I_n$ and $IV$ are impossible.

\begin{table}[t]
\caption{\small Possible singular fibers for the induced morphism of a line $l$ of the second kind of degree 2.}
\label{tab:degree2}
\begin{tabular}{llllll} \toprule
ramification & $F$ & $G$ & $G'$ & $e(G)-e(G')$ & $v_t(l)$ \\
\toprule
\multirow{4}{*}{\small unramified}  & $I_n$      & $2\times I_n$ & $2\times I_{3n}$ & $+4\,n$ & 0 \\
                                    & $I_{3n}$   & $2\times I_{3n}$ & $2\times I_{n}$ & $ - 4\,n$ & 2 \\
                                    & $IV$       & $2\times IV$ & $2\times IV^*$ & $0$ & 2 \\
                                    & $IV^*$     & $2\times IV^*$ & $2\times IV^*$ & $0$ & 2 \\
\midrule
\multirow{5}{*}{\small ramified}    & $I_n^*$    & $I_{2n}$ & $I_{6n}$ & $+4\,n$ & $\leq 1$  \\
                                    & $I_{3n}^*$ & $I_{6n}$ & $I_{2n}$&  $-4\,n$ & $\leq 1$  \\
                                    & $II$       & $IV$ & $IV$ & $0$ & 0  \\
                                    & $II^*$     & $IV^*$ & $IV^*$ & 0 & $\leq 1$  \\
                                    & $IV^*$     & $IV$ & $IV$ & $0$ & $\leq 1$  \\
\bottomrule
\end{tabular}
\end{table}

We can write out Table \ref{tab:degree2}. The two ramified fibers have both Euler number $\geq 2$, so the remaining local contribution is less than or equal than 20. Looking at the possible combinations, one can see that we we get a maximum of 10 lines intersecting $l$.
\endproof

\begin{proposition}
If $l \subset X$ is a line of the second kind of degree~$1$ and singularity $2$, then $v(l) \leq 9$.
\end{proposition}
\proof This is the case when $\alpha$ and $\beta$ have two distinct simple roots in common. Up to projective equivalence, we can suppose that the surface is given by equation~\eqref{eq:surfaceX} with
\[
a_{0130} = a_{0112} = a_{1030} = a_{1021} = a_{1003} = a_{0103} = 0 \quad \text{and} \quad a_{0121} = a_{1012} = 1,
\]
so that the two singular points on $l$ are $P = [0:0:0:1]$ and $Q = [0:0:1:0]$. We have chosen coordinates so that the residual cubic in $\Pi_0: x_0 = 0$ has a double intersection with $l$ at $P$ and the residual cubic in $\Pi_1: x_1 = 0$ has a double intersection with $l$ at $Q$.

One can spell out the conditions for $l$ to be a line of the second kind explicitly; in particular, one finds that $a_{0202} = 0$, so the tangent cone at $P$ splits into two planes $\Pi_0$, $\Pi_2$, whose intersection is different from $l$:
\[
 f_2 = \left( a_{2002} x_{0} + a_{1102} x_{1} + x_{2}\right) x_{0}.
\]

Moreover, we can see, using Bruce-Wall's `recognition principle' \cite[Corollary, p. 246]{brucewall} that the point $P$ must be of type $A_n$ with $n\geq 3$, hence we get a `chain' of $n$ exceptional divisors $\Delta_1,\ldots,\Delta_{n}$, with $\Delta_{i}.\Delta_{i+1} =1$ for $i =1,\ldots n-1$, coming from its minimal resolution (the tangent cone at $P$ splits into two different planes, so $P$ cannot be neither of type $A_1$ nor of types $D_i$, $E_i$; the only case that we must rule out using Bruce-Wall's `recognition principle' is $A_2$). 

Since $l$ is not the intersection of $\Pi_0$ and $\Pi_2$, the general residual cubic of the pencil meets one `extremal' component of the chain of exceptional divisors, say $\Delta_n$, hence if the residual cubic in $\Pi_0$ has $n_0$ components, then the corresponding singular fiber has at least $n_0+2$ components (because it must contain the strict transforms of the $n_0$ components of the residual cubic plus $\Delta_1,\ldots,\Delta_{n-1}$); in particular, it has Euler number $e_0\geq n_0+2\geq 3$. The same applies symmetrically to $Q$: the singular fiber corresponding to the plane $\Pi_1$ has Euler number $e_1 \geq 3$.

Let us denote by $p'$ and $q'$ respectively the number of $p$- and $q$-fibers different from $\Pi_0$ and $\Pi_1$. Since neither $\Pi_0$ nor $\Pi_1$ contributes to the valency of $l$, this cannot be higher than $p'+q'$; in fact, in each $p$-fiber at least two lines must run through the singular points $P$ and $Q$, thus not contributing to the valency of $l$. On the other hand, we must have
\[
 3\, p' + 2\, q' \leq 24 - e_0 - e_1 \leq 18.
\]
Therefore, we infer that $v(l) \leq 9$.
\endproof

\begin{proposition}
If $l \subset X$ is a line of the second kind of degree~$1$ and singularity $1$, then $v(l) \leq 11$.
\end{proposition}

\proof In this case $\alpha$ and $\beta$ have one single double root in common. Suppose the common root is $x_2$, corresponding to the singular point $P = [0:0:0:1]$. Up to projective equivalence, we can choose coordinates so that the residual cubic in $\Pi_0: x_0 = 0$ has triple intersection with $l$ at $P$; hence we can suppose that the surface $X$ is given by equation~\eqref{eq:surfaceX} with 
\[
 a_{0103} = a_{1003} = a_{0112} = a_{1012} = a_{0121} = 0.
\] 
Note that both $a_{0130}$ and $a_{1021}$ must be different from zero, or else the line $l$ would have degree 0.

A necessary condition for $l$ to be of the second kind is $a_{0202} = 0$. The tangent cone at $P$ splits then into two planes: both contain $l$ and one of them is $\Pi_0$; in particular, the point $P$ is not of type $A_1$ and, since the residual cubic in $\Pi_0$ has a singular point in $P$, the corresponding fiber -- which does not contribute to the valency of $l$ -- has Euler number at least~2.

Denote by $p'$ and $q'$ respectively the number of $p$- and by $q$-fibers different from $\Pi_0$. The valency of $l$ is not greater than $p'+q'$ and we have the following bound on the Euler number: 
\[
3\,p' + 2\,q' \leq 24 - 2 = 22
\]
Therefore, $v(l)\leq 11$. 
\endproof

\begin{table}[t]
\caption{\small Known bounds for the valency of a line according to its kind, degree and singularity. Sharp bounds are marked with an asterisk *.}
\label{tab:summary}
 \begin{tabular}{llll} \toprule
 kind & degree & singularity & valency \\
 \midrule
 \multirow{3}{*}{first kind}  & $3$ & 0         & $\leq 18$* \\
                              & $2$ & 1         & $\leq 13$ \\
                              & $1$ & 2 or 1    & $\leq  8$ \\
 \midrule
 \multirow{4}{*}{second kind} & $3$ & 0         & $\leq 20$* \\
                              & $2$ & 1         & $\leq 10$ \\
                              & $1$ & 2         & $\leq 9$  \\
                              & $1$ & 1         & $\leq 11$ \\
 \midrule
 --                           & $0$ & 3, 2 or 1 & $\leq 2$* \\
 \bottomrule
 \end{tabular}
\end{table}

The results of this section are summarized in Table \ref{tab:summary}. Most of the bounds are not known to be sharp; anyway, they are enough for the scope of this article.

\section{Twin lines} \label{sec:twinlines}

We present here a newly discovered configuration of lines, which is also related to the notion of torsion-sections of the Mordell-Weil group. This construction is crucial to the new proof of Segre--Rams--Schütt theorem and to its extension to the K3 quartic case.

\begin{definition}
Let $l\subset X$ a line on a K3 quartic surface. A line $m$ is called an \emph{inflective section} of $l$ if $m$ meets the general residual cubic relative to $l$ in an inflection point.
\end{definition}

\begin{proposition} \label{prop:twinlines}
 Let $X$ be a K3 quartic surface containing two disjoint lines $l$ and $l^*$ of singularity $0$. Then, the following conditions are equivalent:
 \begin{enumerate}[(a)]
  \item There are at least 9 lines $b_1,\ldots,b_9$ meeting $l$ and $l^*$.
  \item There are exactly 10 lines $b_1,\ldots,\,b_{10}$ meeting $l$ and $l^*$.
  \item The line $l^*$ is an inflective section of $l$ and, vice versa, the line $l$ is an inflective section of $l^*$.
  \item The tangents to the general residual cubic $E$ relative to $l$ at the points of intersection of $E$ with $l$ meet in the point of intersection of $E$ with~$l^*$.
  \item The tangents to the general residual cubic $E$ relative to $l^*$ at the points of intersection of $E$ with $l^*$ meet in the point of intersection of $E$ with~$l$.
  \item The quartic $X$ is projectively equivalent to a quartic in the following family $\mathcal A$, where the lines $l$ and $l^*$ are given, respectively, by $x_0 = x_1 = 0$ and $x_2 = x_3 = 0$, and $p_0,\ldots,\,p_3$ are forms of degree 3:
  \begin{equation} \label{eq:familyA}
   \mathcal A : = x_{0} p_{0}(x_{2},x_{3}) + x_{1} p_{1}(x_{2},x_{3}) + x_{2} p_{2} (x_{0}, x_{1}) + x_{3} p_{3} (x_{0}, x_{1})
  \end{equation}
\end{enumerate}
If these conditions are satisfied, then the lines $b_i$ are pairwise disjoint. Moreover, the base change along $l$ induces three 2-torsion sections on $Z$, if one chooses the 0-section to be the one induced by $l^*$, and, vice versa, the base change along $l^*$ also induces three 2-torsion sections on $Z$, if one chooses the 0-section to be the one induced by $l$.
\end{proposition}
\proof Up to coordinate change, we can always suppose that $l$ and $l^*$ are respectively given by $x_0 = x_1 = 0$ and $x_2 = x_3 = 0$.

(a) $\Rightarrow$ (c). The condition of being an inflective section can be computed explicitly and is given by a polynomial of degree 8. The fact that there are at least 9 roots of this polynomial means that it must vanish identically.

(b) $\Rightarrow$ (a) is obvious.

(c) $\Rightarrow$ (f). One computes explicitly the conditions for the lines $l$ and $l^*$ to be inflective sections of each other, and sees that the following coefficients must be equal to zero:
\begin{equation} \label{eq:coeff-twinlines}
    a_{2020},\, a_{2011},\, a_{2002},\, a_{1120},\, a_{1111},\, a_{1102},\, a_{0220},\, a_{0211},\, a_{0202},
\end{equation}
thus obtaining family $\mathcal A$.

(d) $\Rightarrow$ (f). Let $E_t = E_t(x_1,x_2,x_3)$, $t \in \PP^1$, be the residual cubic relative to $l$, obtained by substituting $x_0 = t x_1$ in the equation of $X$. Consider the polynomial $\partial E_t/\partial x_1$ restricted on the line $l:x_0 = x_1 = 0$: this is a polynomial of degree 2 in $(x_2,x_3)$ with polynomials of degree 2 in $t$ as coefficients. Since generically it must have three distinct roots, namely the points of intersection of $E_t$ with $l$, it must be the zero polynomial; hence, the coefficients must be the zero polynomial in $t$. The result is that the same coefficients listed in \eqref{eq:coeff-twinlines} must be equal to zero.

(e) $\Rightarrow$ (f) is proven analogously.

(f) $\Rightarrow$ (c), (d), (e) is immediate.

(f) $\Rightarrow$ (b) can be proven explicitly by considering the discriminant of the fibration induced by one of the two lines.
\endproof

\begin{definition}
    If $l$ and $l^*$ satisfy one of the equivalent conditions of Proposition \ref{prop:twinlines}, we say that $l$ and $l^*$ are \emph{twin lines}.
\end{definition}

\begin{remark}
    The family $\mathcal A$ has dimension 8; in fact, knowing that there are 10 disjoint lines meeting both $l$ and $l^*$, we can assume -- up to projective equivalence -- that two of them are given respectively by $x_1 = x_2 = 0$ and $x_0 = x_3 = 0$; we are left with 12 parameters, 4 of which can be normalized to 1. Indeed, the lattice generated by the twelve lines and the hyperplane section has rank 12, as expected.
\end{remark}

\begin{remark} \label{rmk:tau}
The explicit parametrization \eqref{eq:familyA} of family $\mathcal A$ shows the existence of a (non-symplectic) automorphism $\tau:X \rightarrow X$ of degree 2, given by
\[
    \tau: [x_0:x_1:x_2:x_3] \mapsto [-x_0:-x_1:x_2:x_3].
\]
This automorphism fixes $l$ and $l^*$ pointwise; it also respects their fibers as sets.
\end{remark}

\begin{corollary} \label{cor:p-fiber_twin}
    Let $l \subset X$ a line of singularity 0 admitting a twin $l^*$. If $l$ has a $p$-fiber, then this fiber is ramified.
\end{corollary}
\proof Suppose $X$, $l$ and $l^*$ are given as in family $\mathcal A$. Exactly one of the lines in the $p$-fiber of $l$ must meet the line $l^*$: let us call it $m_0$ and the other two $m_1$ and $m_2$. The points of intersection of $m_0$ with $l$ and $l^*$ are fixed by the automorphism $\tau$ (see Remark \ref{rmk:tau}); hence, $m_0$ is mapped to itself.

Necessarily, the point $P$ of intersection of $m_1$ and $m_2$ is also fixed by $\tau$. Since $P$ does not lie on $l^*$, one of its last two coordinates must be different from zero; this implies that its first two coordinates must be zero; therefore, it must lie on $l$ and ramification must occur (of index 3 or 2, according to whether the lines $m_i$ meet at the same point or not). \endproof

\begin{lemma}[Degtyarev--Itenberg--Sertöz] \label{lemma:(4,6)}
 If $l\subset X$ is a line of singularity~$0$ inducing a fibration of type $(p,q) = (4,6)$, then $X$ is smooth and the line $l$ has a twin $l^*$.
\end{lemma}
\proof[Sketch of proof] There cannot be singular points outside $l$, otherwise the Euler number of $X$ would exceed 24; since $l$ has no singular points, the surface is smooth. The surface $X$ is therefore a K3 surface. Let us call $m_{i,j}$, $i=1,\ldots,4$, $j=1,2,3$, the lines in the $p$-fibers and $n_k$, $k=1,\ldots,6$, the lines in the $q$-fibers. 

Suppose first that the base field has characteristic $0$. The lattice $L$ generated by the lines and the hyperplane section must admit an embedding into the K3 lattice $\Lambda = U^3 \oplus E_8(-1)^2$. By results of Nikulin \cite{nikulin}, this embedding cannot be primitive, due to the $3$-elementary part of the discriminant group of $L$. A careful analysis of the admissible isotropic vectors reveals that -- up to symmetry -- the following class must also be contained in the Picard lattice of $X$:
\[
\omega : = \frac 13 \left(l + \sum_{i=1}^4 (m_{i,1} + m_{i,2}) - \sum_{k=1}^6 n_k \right).
\]
One can check that this is exactly the class of the sought line $l^*$.

If the base field has positive characteristic $p>3$, one has to distinguish two cases. 
\begin{itemize}
\item If the surface is not Shioda-supersingular, then one can lift it -- together with the whole Picard group -- to characteristic 0 (see, for instance, Esnault--Srinivas \cite[p. 839]{esnault-srinivas}), so that one can apply the same arguments.
\item If the surface is Shioda-supersingular, then the lattice $L$ must embed in a $p$-elementary lattice. Since $p>3$, one obtains the same condition on the $3$-elementary part of the discriminant group of $L$ which prevents it from embedding primitively. Again, one concludes that $\omega$ must be contained in the Picard lattice. \qedhere
\end{itemize}
\endproof

\begin{remark}
Our interest in lines of type $(4,6)$ was motivated by the fact that Schur's quartic contains 16 lines of type $(6,0)$ of the second kind and 48 lines of type $(4,6)$ of the first kind. We point out a mistake in Rams and Schütt's article \cite{64lines}: Proposition 7.1, which claims that in a quartic containing 64 lines all lines are of type $(6,0)$, is false; the flaw lies in the proof of Lemma 7.3 [\emph{ibidem}].
\end{remark}

\section{Lines of high valency} \label{sec:highvalency}

Table \ref{tab:summary} and Proposition \ref{prop:degree3} prompt us to give a closer look at surfaces containing a line of degree 3 admitting only two points of ramification, which are the only ones that might have valency greater than 19. We can parametrize such surfaces in the same way as Rams and Schütt did \cite[Lemma 4.5]{64lines}.

\begin{lemma} \label{lemma:familyZ}
 Let $X \subset \PP^3$ be a K3 quartic surface containing a line $l$ of the second kind which induces a morphism $\pi: \hat l \rightarrow \PP^1$ of degree 3 ramifying over two points with ramification index 3. Then, $X$ is projectively equivalent to a quartic in the family
 \begin{equation} \label{eq:parametrizationZ}
  \mathcal Z: x_0 x_3^3 + x_1 x_2^3 + x_2 x_3 q_2(x_0, x_1) + q_4(x_0,x_1) = 0,
 \end{equation}
 where $q_i \in k[x_0, x_1]$ are homogeneous polynomials of degree $i$ ($i = 2,4$).
\end{lemma}
\proof Knowing that there are no singular points on the line $l$, the proof can be copied word by word from \cite[Lemma 4.5]{64lines}. In the proof one uses the fact that the characteristic of the ground field is different from 3. \endproof

Notably, Schur's quartic, given by the equation 
\begin{equation} \label{eq:Schur}
 x_0^4 - x_0 x_3^3 = x_1^4 - x_1 x_2^3,
\end{equation}
is projectively equivalent to a member of the family $\mathcal Z$. It is well known that it contains exactly 64 lines over any algebraically closed field of characteristic different from 2 and 3 (see, for example, Boissière-Sarti~\cite{boissiere-sarti}).

\begin{remark} \label{rem:sigma}
    The parametrization given by equation \eqref{eq:parametrizationZ} reveals that there exists a (symplectic) automorphism $\sigma:X \rightarrow X$ of order 3 which is given by
    \[
    \sigma: [x_0:x_1:x_2:x_3] \mapsto [x_0: x_1: \zeta x_2: \zeta^2 x_3],
    \]
    with $\zeta$ a primitive third root of unity. In what follows we will refer to this automorphisms as `the' automorphism of order 3 induced by $l$. Note that $\sigma$ permutes the components of the $p$-fibers of $l$.
\end{remark}

\begin{remark}
    It is worth noticing that if a K3 quartic $X$ contains a line with valency greater than 18, then $X$ is projectively equivalent to a member of the family $\mathcal Z$. This follows from Table \ref{tab:summary} and Proposition \ref{prop:degree3}, and parallels \cite[Proposition 1.1]{64lines}. On the other hand, if one considers the \emph{extended} valency of $l$, then there are more configurations where this can be greater than 18; examples of such phenomena are illustrated in Section \ref{sec:examples}.
\end{remark}

The following proposition is a generalization of \cite[Lemma 6.2]{64lines} to the K3 quartic case.
\begin{proposition} \label{prop:secondkindimplies16}
    Let $l,\,m$ be two intersecting lines on a K3 surface $X$. Suppose that both induced morphisms $\hat l \rightarrow \PP^1$ and $\hat m \rightarrow \PP^1$ have degree 3 and exactly two points of ramification, and that both $l$ and $m$ are of the second kind. Then, $X$ is projectively equivalent to Schur's quartic.
\end{proposition}
\proof Let $P$ be the point of intersection of $l$ and $m$; $Q$ one of the two ramification points on $l$, corresponding to the plane $\Pi \supset l$; $R$ one of the ramification points of $m$, corresponding to the plane $\Sigma \supset m$; $S$ one of the points of intersection of the line $\Pi \cap \Sigma$ with $X$ different from $P$.

Up to projective equivalence, we can suppose that $P$, $Q$, $R$ and $S$ are respectively the points $[0:0:1:0]$, $[0:0:0:1]$, $[0:1:0:0]$ and $[1:0:0:0]$. Thus, the line $l$ is given by $x_0 = x_1 = 0$ and the line $m$ by $x_0 = x_3 = 0$.

This amounts to setting the following coefficients equal to zero in equation~\eqref{eq:surfaceX}:
\[
    a_{0400},\, a_{0310},\, a_{0220},\, a_{0130},\, a_{1003},\, a_{1012},\, a_{1021},\, a_{1300},\, a_{1210},\, a_{1120},\, a_{4000}.
\]
Furthermore, since $l$ and $m$ do not contain singular points, the following coefficients must be different from zero and we can set them to 1:
\[
    a_{0103},\,a_{0301},\, a_{1030}.
\]

Recall that a necessary condition for a cubic polynomial 
  \[
   p(t) = a t^3 + b t^2  + c t  + d
  \]
to have a triple root is
  \[
   b^2 - 3 \, a c = 0.
  \]
  Therefore, in order for $m$ and $l$ to have exactly two points of ramification, one sees that the following equations must be satisfied:
  \[
  3\,a_{0121} = a_{0211}^2 = a_{0112}^2.
  \]
Spelling out the conditions for $l$ and $m$ to be of the second kind, one sees that these coefficients must be actually zero. Indeed, one obtains a surface which is immediately seen to be projectively equivalent to Schur's quartic.
\endproof

\begin{corollary} \label{cor:intersecting-lines-of-high-valency}
    A K3 quartic surface cannot contain two intersecting lines of valency greater than 18.
\end{corollary}
\proof On account of Table \ref{tab:summary} and Proposition \ref{prop:degree3}, the two lines would be as in Proposition \ref{prop:secondkindimplies16}, but all lines on Schur's quartic have valency 18.
\endproof

\begin{proposition} \label{prop:p-fiber-singular-point}
    Let $X$ be a K3 quartic surface with a singular point $P$. Suppose a plane $\Pi$ containing $P$ splits $X$ into four (not necessarily distinct) lines. Then, the surface $X$ contains at most 62 lines.
\end{proposition}
\proof One can list the possible configurations of lines and singular points that can lie on the plane $\Pi$; then, one uses the bound of 8 lines through each singular point and the bounds given by Table \ref{tab:summary}. In almost every configurations where there is a line $l$ of singularity~0 one can assume that $v(l)\leq 18$; in fact, if this is not the case then we can infer from the presence of the automorphism $\sigma$ that either the other three lines meet at the same singular point, or they meet at three distinct singular points, or they coincide.

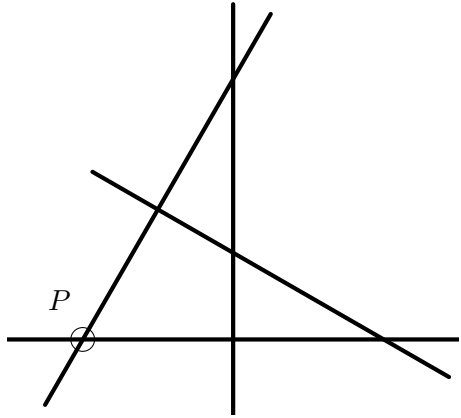
\begin{figure}
\definecolor{uququq}{rgb}{0,0,0}
\begin{tikzpicture}[line cap=round,line join=round,>=triangle 45,x=1.0cm,y=1.0cm] 
\clip(-1,-1) rectangle (5,5);
\draw [line width=1.6pt] (-0.5,-0.87)-- (2.5,4.33);
\draw [line width=1.6pt] (0.13,2.23)-- (4.87,-0.5);
\draw [line width=1.6pt] (-1,0)-- (5,0);
\draw [line width=1.6pt] (2,4.46)-- (2,-1);
\draw [color=uququq] (0,0) circle (4.5pt);
\draw[color=black] (-0.31,0.52) node {$P$};
\end{tikzpicture}
\caption{\small A plane splitting $X$ into four lines with a singular point $P$.}
\label{fig:p-fiber-singular-point}
\end{figure}

The case admitting the worst bound is when there is one singular point $P$ through two lines, as pictured in Figure \ref{fig:p-fiber-singular-point}. In fact, there might be up to 8 lines through $P$,  while the valency of the two lines through $P$ is not greater than 13 and the valency of the other two lines is not greater than 18, yielding a total bound of
\[
    (8-2) + 2\cdot(13 - 2) + 2\cdot(18-3) + 4 = 62. \qedhere
\]

\endproof

\section{Triangle-free quartic surfaces} \label{sec:trianglefree}
Most of the ideas contained in this section are due to A. Degtyarev.

\begin{definition}
    Let $X$ be a K3 quartic surface with minimal desingularization $Z$. The \emph{line graph} of $X$ is the dual graph of the lines on $Z$, i.e. it is the graph whose vertex set is the set of lines on $X$ such that two vertices $l$, $m$ are connected by an edge if and only if $\hat l \cdot \hat m = 1$, where $\hat l$, $\hat m$ are the strict transforms of the lines $l$ and $m$ in $Z$.
\end{definition}
The line graph of a K3 quartic surface is a graph without loops or multiple edges.

\begin{definition}
    A K3 quartic surface $X$ is called \emph{triangle-free} if its line graph contains no triangles, i.e. cycles of length 3.
\end{definition}

In other words, a K3 quartic surface is triangle-free if there are no triples of pairwise intersecting lines on its minimal desingularization $Z$; equivalently, a K3 quartic surface $X$ is triangle-free of there are no triples of lines on $X$ meeting pairwise at smooth points. The next definition has an analogous geometric interpretation.

\begin{definition}
    A K3 quartic surface $X$ is called \emph{quadrangle-free} if it is triangle-free and if its line graph contains no quadrangles, i.e. cycles of length 4.
\end{definition}

Recall that a graph induces a symmetric bilinear form on the lattice generated by its vertices, in the following way:
\begin{align*}
 v^2 = v \cdot v & := -2 + 2\cdot\#\{\text{loops around $v$}\} \\
 v \cdot w & := \#\{\text{edges joining $v$ and $w$}\}
\end{align*}
Note that the symmetric form on the line graph of a K3 quartic coincides with the intersection form on the lines contained in its minimal desingularization.

\begin{definition} 
  A connected graph is called \emph{elliptic} if its associated form is negative definite; \emph{parabolic} if its associated form is negative semidefinite, with kernel of dimension 1. In other words, elliptic graphs are Dynkin diagrams and parabolic graphs are extended Dynkin diagrams.
\end{definition}

In what follows, by `subgraph' we will always mean an `induced subgraph'. We will denote by $|G|$ the cardinality of the set of vertices of a graph $G$. The Milnor number $\mu(G)$ of a graph $G$ is the rank of its associated form.

Let now $\Gamma$ be the line graph of a K3 quartic surface $X$. Given a subgraph $G\subset \Gamma$, the \emph{span} of $G$ will be the subgraph
\[
 \Span G = G \cup \{m \in \Gamma: m\cdot l = 1 \text{ for some $l \in G$}\};
\]
the \emph{valency} of $G$ will be 
\[
 v(G) := |(\Span G) \smallsetminus G|.
\]
Note that this definition extends naturally the notion of `valency of a line'.

Two subgraphs $G,\,G'$ of $\Gamma$ are said to be \emph{disjoint} if $\Span G \cap G' = \varnothing$. 

\begin{proposition} \label{prop:alex}
    Let $\Gamma$ the line graph of a K3 quartic surface $X$. If $\Gamma$ contains a parabolic subgraph $D$, then 
    \[
      |\Gamma| \leq v(D) + 24.
    \]
\end{proposition}
\proof A parabolic subgraph induces an elliptic fibration \cite[§3, Theorem 1]{PSS}. The vertices in $D \cup (\Gamma \smallsetminus \Span D)$ are fiber components of this fibration; hence, on account of the Euler number, they cannot be more than 24 in number. \endproof

We now set
\[
 \delta := \begin{cases} 20 & \mbox{if } \characteristic k = 0 \\
			 22 & \mbox{if } \characteristic k > 0. \end{cases}
\]
The number $\delta$ is a well-known bound for the rank of the Picard lattice of a K3 surface. Since the signature of the Picard lattice is $(1,\rho-1)$, the Milnor number of any negative semidefinite subgraph of $\Gamma$ cannot be greater than $\delta -1$. In particular, since $\Gamma$ has neither loops nor multiple edges, it can only contain the following parabolic subgraphs: $\tilde A_2,\ldots,\,\tilde A_{\delta - 1}$, $\tilde D_n,\ldots,\,\tilde D_{\delta-1}$, $\tilde E_6$, $\tilde E_7$, $\tilde E_8$.

\begin{lemma}
 If $\Gamma$ does not contain any parabolic subgraph, then $|\Gamma| \leq \delta - 1$.
\end{lemma}
\proof
The associated form of $\Gamma$ must be negative definite; hence, $\Gamma$ is the disjoint union of elliptic graphs and its Milnor number is equal to $|\Gamma|$. 
\endproof

\begin{lemma} \label{lemma:v-leq3}
 If $v(l) \leq 3$ for every line $l\subset X$, then $|\Gamma| \leq \delta + 26.$
\end{lemma}
\proof By the previous lemma, we can assume that there is a parabolic subgraph $D\subset \Gamma$. Under the hypothesis $v(l) \leq 3$ for every vertex $l \in \Gamma$, we deduce that $v(\tilde A_n) \leq n+1$, $v(\tilde D_n) \leq n+3$, $v(\tilde E_n) \leq n+3$; hence, by virtue of Proposition \ref{prop:alex}, we obtain
\[
 |\Gamma| \leq v(D) + 24 \leq (n+3) +24 \leq (\delta -1 +3)+24 = \delta +26. \qedhere
\]
\endproof

\begin{proposition} \label{prop:quadrangle-free}
    A quadrangle-free K3 quartic surface contains at most 54 lines (51 over a field of characteristic 0).
\end{proposition}
\proof
Let $l$ be a vertex of maximal valency. On account of Lemma~\ref{lemma:v-leq3}, we can assume that the valency $w$ of $l$ is at least 4. Let $m_1,\ldots,\,m_4$ be four vertices adjacent to $l$ and suppose that $m_i$ has valency $v_i$ ($i = 1,\ldots,4$).

Since the surface is quadrangle-free, all vertices adjacent to $l$ are disjoint from the vertices adjacent to $m_i$; moreover, a vertex adjacent to $m_i$ can be joined to at most one vertex adjacent to $m_j$, for $i\neq j$. Hence, we can assume that there are $a$ lines meeting $l$ different from $m_i$, $m_j$; $b$ lines meeting $m_i$ not intersecting any line meeting $m_j$; $c$ lines meeting $m_j$ not intersecting any line meeting $m_i$; $d$ pairs of lines forming a pentagon with $m_i$, $m_j$ and $l$. Note that
\[
    s := w + v_i + v_j = a + b + c + 2\,d + 4.
\]

A simple computer-aided computation -- which amounts to constructing all possible intersection matrices of $\Span \{l,m_i,m_j\}$ for values of $a,\,b,\,c,\,d$ such that $s = \delta + 2$ and computing their ranks -- shows that any configuration with $s \geq \delta + 2$ gives rise to a lattice of rank greater than $\delta$. Thus, we can assume that
\begin{equation} \label{eq:w1}
    w + v_i + v_j  \leq \delta + 1 \quad \text{for all $i \neq j$}.
\end{equation}
Taking the sum of \eqref{eq:w1} with $(i,j) = (1,2),\,(3,4)$, one finds that
\begin{equation} \label{eq:w}
    w + \sum_{i=1}^4 v_i \leq 2\,(\delta + 1) - w.
\end{equation}

Hence, if $w \leq \lfloor (\delta + 1)/3 \rfloor$, by the maximality condition one has
\[
    w + \sum_{i=1}^4 v_i \leq 5\, w \leq \frac 5 3 (\delta + 1);
\]
on the other hand, if $w \geq \lceil (\delta + 1)/3 \rceil$, then one obtains the same relation from inequality~\eqref{eq:w}. Applying Proposition~\ref{prop:alex} to the $\tilde D_4$-subgraph formed by $l$, $m_1,\ldots,\,m_4$, we find that
\[
    |\Gamma| \leq (w - 4) + \sum_{i=1}^4 (v_i - 1) + 24 \leq 16 + \frac 5 3 (\delta + 1),
\]
which yields the claim.
\endproof

\begin{proposition} \label{prop:triangle-free}
    A triangle-free K3 quartic surface contains at most 64 lines.
\end{proposition}
\proof By virtue of Proposition \ref{prop:quadrangle-free}, we can assume that $\Gamma$ contains a quadrangle $D$, formed by the lines $l_i$, $i=1,\ldots,4$. If one of the $l_i$ had a $p$-fiber, then the corresponding plane should contain at least one singular point, because the surface is triangle-free; hence, we could apply Proposition~\ref{prop:p-fiber-singular-point}. Therefore, we can suppose that the $l_i$'s have no $p$-fibers; thus, by virtue of Lemma~\ref{lemma:no-p-fibers}, they must have valency at most $12$.

Applying Proposition \ref{prop:alex} to the $\tilde A_3$-subgraph $D$, we infer that
\[
    |\Gamma| \leq 4\cdot(12 - 2) + 24 = 64. \qedhere
\]
\endproof

\begin{remark}
    The bounds presented in this section are most probably not sharp. We were informed by A. Degtyarev that he has found a triangle-free smooth surface over $\CC$ with 37 lines.
\end{remark}

\section{Proof of main theorem} \label{sec:prooftheorem}

We are now able to present a proof of Theorem \ref{thm:maintheorem}, which we restate here in a more concise form.

\begin{theorem}
    A K3 quartic surface $X$ contains at most 64 lines.
\end{theorem}
\proof Assume that $X$ contains more than 64 lines. By virtue of Proposition~\ref{prop:triangle-free}, we can suppose that $X$ is not triangle-free; hence, there are three lines $l_1$, $l_2$ and $l_3$ on $X$ that meet at smooth points of $X$. These three lines are all contained in a plane $\Pi$; let $l$ be the fourth line on this plane. On account of Proposition~\ref{prop:p-fiber-singular-point}, we can assume that all points on $\Pi$ are smooth.

If all four lines on $\Pi$ have valency less or equal than $18$, then the total number of lines lying on $X$ can be at most
\[
    4 + 4 \cdot (18 - 3) = 64.
\]
Hence, we can assume that one of the lines, say $l$, has valency $19$ or $20$. By virtue of Table \ref{tab:summary} and Proposition \ref{prop:degree3}, the line $l$ has a ramified $q$-fiber. Let us call $m$ and $C$ respectively the line and the conic in the $q$-fiber, and $P$ the point of intersection of $m$ and $C$ not lying on~$l$. Note that $P$ may or may not be a singular point of~$X$.

\begin{claim} \label{claim:vm}
 $v(m) \leq 10$.
\end{claim}
\proof[Proof of the claim] If $P$ is not singular, then the line $m$ is of the first kind, since the point $P$ is certainly not an inflection point of the corresponding residual cubic, whence $v(m) \leq 18$. On account of the automorphism $\sigma$ of order 3 induced by $l$ (see Remark~\ref{rem:sigma}), we know that the valency of $m$ has the form
\[
v(m) = 1 + 3\,a,
\]
for some integer $a \geq 0$; thus, $v(m)\leq 16$.  On the other hand, $m$ has no $p$-fibers: in fact, since $v(m) \leq 16$, if $m$ had a $p$-fiber (with no singular points -- see Proposition~\ref{prop:p-fiber-singular-point}), then at least one line $n$ in the $p$-fiber should have valency 19 or 20; the automorphism $\sigma'$ induced by $n$ would force the other two residual lines to have the same valency as $m$, that is to say, not greater than 16: all in all, the lines on $X$ would be less than 64, which is absurd. Hence, by Lemma \ref{lemma:no-p-fibers}), it follows that $v(m) \leq 12$, i.e. $a \leq 3$.

On the other hand, if $P$ is a singular point, it follows from Proposition~\ref{prop:p-fiber-singular-point} that $m$ cannot have $p$-fibers, since we are assuming that $X$ contains more than 64 lines.
Given that in this case, too, the valency of $m$ has the form $1+3\,a$, we can conclude as before.
\endproof

\begin{claim}
 The lines $l_i$ ($i = 1,\,2,\,3$) are of type $(4,6)$.
\end{claim}

\proof[Proof of the claim] Because of the presence of the automorphism $\sigma$, the lines in question have the same valency $v$ and are of the same type. On the one hand, $v$ cannot be greater than 18, by Corollary \ref{cor:intersecting-lines-of-high-valency}; on the other hand, if $v \leq 17$, then the total number of lines on $X$ would be at most
\[
    4 + 3 \cdot (v - 3) + (v(l) - 3) \leq 4 + 3 \cdot (17 - 3) + (20 - 3) = 63.
\]
Therefore, $v$ is exactly 18, whence the $l_i$'s must have type $(p,q)=(6,0)$, $(5,3)$ or $(4,6)$; in fact, $p \leq 3$ is not possible, by a simple Euler number argument.

Observe now that if $P$ is singular, then the plane containing $l_i$ and $P$ cannot be a $p$-fiber for $l_i$, by Proposition~\ref{prop:p-fiber-singular-point}. Regardless of $P$ being smooth or singular, it follows that all $p$-fibers of $l_1$, $l_2$ and $l_3$ contain a line meeting $m$ in a point different from $P$. Since by Claim~\ref{claim:vm} the valency of $m$ is at most 10, the $l_i$'s can only have type $(4,6)$.
\endproof

By virtue of Corollary~\ref{cor:p-fiber_twin}, the plane $\Pi$ is a ramified fiber of ramification index 2 for the lines $l_i$; hence, the three lines must meet in a point. In particular, the plane $\Pi$ is a fiber of type $IV$ for the line $l$. Recall that, by Proposition~\ref{prop:degree3}, the line $l$ has 6 $p$-fibers; considering that we could repeat the same argument for any $p$-fiber of $l$, we deduce that $l$ has 6 fibers of type $IV$. But since the line $l$ has also at least one $q$-fiber, we deduce that the Euler number of the minimal desingularization of $X$ must be at least
\[
 6 \cdot 4 + 2 = 26,
\]
which is impossible. \endproof

\section{Non-smooth surfaces with many lines} \label{sec:examples}

In Section \ref{sec:valency} we have devoted ourselves to finding bounds for the valency of a line $l$ contained in a K3 quartic surface $X$; Table \ref{tab:summary} shows that $v(l)$ can never be greater than 20 and only in very special configurations it can be greater than 18. Another natural question is the following: what is the maximal \emph{extended} valency that $l$ can have? A rough answer to this question is that, since there cannot be more than 8 lines through a singular point by Lemma \ref{lemma:linesthroughsingularpoint},
\[
\tilde v(l) \leq v(l) + 7\,s,
\]
where $s$ denotes the singularity of the line $l$.

A much more careful analysis reveals that also
\[
    \tilde v(l) \leq 20,
\]
for any line $l\subset X$. Apart from the obvious case of $l$ as in Proposition \ref{prop:degree3}, where $\tilde v(l) = v(l)$ can be greater than 18, there are two other configurations in which $\tilde v(l) >18$. We do not present a proof of these assertions, as we were only interested in giving a (as sleek as possible) proof of Theorem \ref{thm:maintheorem}, but we provide examples of surfaces which display these new behaviors. These examples were found by examining lines of the first kind and lines of degree~0 more closely, for example taking advantage of the restrictions of Lemma~\ref{lemma:linesthroughsingularpoint}.

\begin{example}
 Consider the surface defined by
\begin{align*}
\begin{split}
   & 3 \, x_{0}^{4} - 9 \, x_{0}^{3} x_{1} + 6 \, x_{0}^{2} x_{1}^{2} - 12 \, x_{0} x_{1}^{3} + 8 \, x_{1}^{4} - 9 \, x_{1}^{3} x_{2} +\\ 
   & - 27 \, x_{0}^{2} x_{2}^{2} - 27 \, x_{1}^{2} x_{2}^{2} - 27 \, x_{1} x_{2}^{3} - 27 \, x_{1}^{2} x_{3}^{2} - 27 \, x_{0} x_{2} x_{3}^{2} = 0.
\end{split}
\end{align*}
It has one singular point $P$ of type $A_1$. The line $l$ given by $x_0 = x_1 = 0$ is a line of the first kind of singularity 1 and degree 2; it has valency 12 and extended valency 19.

The fibration induced on the minimal desingularization $Z$ by $l$ has six fibers $F_i$ of type $I_3$ and one fiber $G$ of type $I_2$. There are no other lines on the surface; hence, the surface contains exactly 20 lines. This holds true over any field of characteristic $p\neq 2,3$ such that this fibration does not degenerate, and one can check it in the following way.

The fibers $F_i$ come from residual cubics composed of three lines: we call $m_{i,0}$ the line passing through $P$, and $m_{i,1},m_{i,2}$ the other two lines ($i = 0,\ldots,5$). We call $n$ the line in the residual cubic corresponding to $G$, which passes through $P$. One can check explicitly that the lines $m_{i,0}$ and $n$ do not meet other lines. The intersection matrix of the strict transforms of all these lines on $Z$ and of the exceptional divisor resulting from the blowup of $P$ has signature $(1,14)$. If a section $s$ existed, then it would meet exactly one line between $m_{i,1},m_{i_2}$, for $i=0,\ldots,5$, and no other line; up to symmetry, we can suppose that $s$ would meet $m_{i,1}$ for $i=0,\ldots,5$. However, the resulting intersection matrix would have signature $(1,16)$, which is impossible, since adding one divisor the signature must either stay the same or become $(1,15)$.
\end{example}

\begin{example}[due to González Alonso and Rams] \label{ex:Gonz-Rams}
 The surface over $\CC$ given by
 \[
  x_{0}^{4} + x_{0} x_{2}^{3} + x_{1}^{2} x_{2} x_{3} + x_{0} x_{3}^{3} = 0
 \]
 has one singular point of type $A_3$ and 3 singular points of type $A_1$. The line given by $x_0 = x_1 = 0$ has singularity 3, valency 2 and extended valency 20.

 This surface contains exactly 39 lines. To our knowledge, this is so far the example of an explicit non-smooth K3 quartic surface with most lines over a field of characteristic zero. González Alonso and Rams came to this example by checking all Delsarte surfaces in Heijne's list \cite{heijne}.

 By a careful inspection of the fibrations induced by the lines lying in the plane $x_0 = 0$, one can conclude that there exists no prime $p$ such that the reduction of this surface modulo $p$ contains more than 39 lines.
\end{example}

\begin{example} \label{ex:implicit}
A non-smooth complex K3 quartic surface with 40 lines exist and has been found with Degtyarev-Itenberg-Sertöz' lattice-theoretical methods \cite{AIS}. It contains one singular point of type $A_1$. An explicit equation of the surface is not known.
\end{example}

The bound expressed by Theorem~\ref{thm:maintheorem} is sharp, since Schur's quartic \eqref{eq:Schur} -- which is smooth -- contains exactly 64 lines. It is still an open question what the maximum number of lines on non-smooth K3 quartic surfaces is. In other words, given a quartic surface $X\subset \PP^3$ with at most isolated rational double points, but with at least one singular point, what is the maximum number of lines that can lie on $X$?

A part from the Examples \ref{ex:Gonz-Rams} and \ref{ex:implicit} with 39 and 40 lines over $\CC$, we list here some notable surfaces with many lines defined over fields of positive characteristic, with 42, 45 and 48 lines. The following examples were found either by inspecting the family $\mathcal{Z}$ \eqref{eq:parametrizationZ} or by imposing a lot of symmetries on the surface. It is worth mentioning here that González Alonso and Rams proved that a complex non-ruled quartic surface with worse singularities than isolated rational double points (i.e. a complex non-K3 quartic surfaces not containing an infinite number of lines) can contain at most 48 lines and their best example (due to Rohn) contains 31 lines and a triple point \cite{GA-Rams}.

\begin{example}
 The surface over $\CC$ given by 
 \[
  x_{0}^{2} x_{1}^{2} + x_{1} x_{2}^{3} - x_{0}^{2} x_{2} x_{3} - x_{0} x_{1} x_{2} x_{3} - x_{1}^{2} x_{2} x_{3} + x_{0} x_{3}^{3} = 0
 \]
 belongs to the family $\mathcal Z$, has 5 singular points of type $A_1$ and contains exactly 33 lines. It has Picard number 20. 
 Its reduction modulo 5 contains 42 lines and has Picard number 22.
\end{example}

\begin{example}
 The surface over $\CC$ given by 
 \[
  x_{0}^{3} x_{1} - 2 \, x_{0}^{2} x_{1}^{2} + x_{0} x_{1}^{3} + x_{1} x_{2}^{3} + x_{0}^{2} x_{2} x_{3} - x_{0} x_{1} x_{2} x_{3} + x_{1}^{2} x_{2} x_{3} + x_{0} x_{3}^{3} = 0
 \]
 belongs to the family $\mathcal Z$, has 1 singular point of type $A_1$ and contains exactly 36 lines. It has Picard number 20. 
 Its reduction modulo 11 contains 45 lines and has Picard number 22.
\end{example}
\begin{example} \label{ex:48lines}
 The surface over $\CC$ given by
 \[
  x_{0}^{2} x_{1} x_{2} + x_{1}^{2} x_{2}^{2} + x_{0} x_{1}^{2} x_{3} + x_{0} x_{2}^{2} x_{3} + x_{0}^{2} x_{3}^{2} + x_{1} x_{2} x_{3}^{2} = 0
 \]
 has 4 singular points of type $A_1$ and contains exactly 36 lines.
 Its reduction modulo 5 contains exactly 48 lines. To our knowledge, this is so far the example of a non-smooth K3 quartic surface with most lines over a field of positive characteristic $p\neq 2,3$.
\end{example}

\section*{Aknowledgements}
First of all, I wish to warmly thank my supervisor Matthias Schütt and S\l awek Rams for suggesting the problem and paving the way to solve it. Many new ideas are due to the fruitful discussions with Alex Degtyarev, who made my stay in Ankara extremely pleasant. Thanks to Miguel Ángel Marco Buzunáriz for introducing me to SageMath, to Roberto Laface for carefully proofreading the draft, and to Víctor González Alonso and Simon Brandhorst for their invaluable mathematical help.


\begin{thebibliography}{30}
 \bibitem{boissiere-sarti} S. Boissière, A. Sarti. \emph{Counting lines on surfaces}. Ann. Scuola Norm. Sup. Pisa, Cl. Sci. \textbf{5} (2007), 39--52.
 \bibitem{brucewall} J. W. Bruce, C. T. C. Wall. \emph{On the classification of cubic surfaces}. J. London Math. Soc. \textbf{19}(2) (1979), 245--256.
 \bibitem{AIS} A. Degtyarev, I. Itenberg, A. S. Sertöz. \emph{Lines on quartic surfaces}. Unpublished.
 \bibitem{esnault-srinivas} H. Esnault, V. Srinivas. \emph{Algebraic versus topological entropy for surfaces over finite fields}. Osaka J. Math.
\textbf{50}(3) (2013), 827--846.
 \bibitem{GA-Rams} V. González Alonso, S. Rams. \emph{Counting lines on quartic surfaces}. Preprint available at arXiv:1505.02018 (May 2015).
 \bibitem{kollar} J. Kollár. \emph{Szémeredi--Trotter-type theorems in dimension 3}. Adv. Math. \textbf{271} (2015), 30--61. 
 \bibitem{heijne} B. Heijne. \emph{Picard numbers of complex Delsarte surfaces with only isolated ADE-singularities.} Preprint available at arXiv:1212.5006v4 (May 2015).
 \bibitem{liedtke} C. Liedtke. \emph{Algebraic Surfaces in Positive Characteristic}. Preprint available at arXiv:0912.4291v4 (May 2015).
 \bibitem{miranda} R. Miranda. \emph{The basic theory of elliptic surfaces}. ETS Editrice Pisa (1989).
 \bibitem{nikulin} V. V. Nikulin. \emph{Integral symmetric bilinear forms and some of their applications}. Math. USSR Izvestija \textbf{14}(1) (1980), 103--167.
 \bibitem{PSS} I. Pjateckij-Šapiro, I. Šafarevič. \emph{A Torelli theorem for algebraic surfaces of type K3}. Izv. Akad. Nauk SSSR \textbf{35} (1971), 530--572.
 \bibitem{64lines} S. Rams, M. Schütt. \emph{64 lines on smooth quartic surfaces}. Math. Ann. \textbf{362}(1) (2015), 679--698.
 \bibitem{onquartics} S. Rams, M. Schütt. \emph{On quartics with lines of the second kind}.  Adv. Geom. \textbf{14} (2014), 735--756. 
 \bibitem{char3} S. Rams, M. Schütt. \emph{112 lines on smooth quartic surfaces (characteristic 3)}. Preprint available at arXiv:1409.7485 (May 2015). To appear in Quart. J. Math.
 \bibitem{saint-donat} B. Saint-Donat. \emph{Projective models of K3 surfaces}. Amer. J. Math. \textbf{96} (1974), 602--639.
 \bibitem{segre} B. Segre. \emph{The maximum number of lines lying on a quartic surface}, Quart. J. Math., Oxford Ser. \textbf{14} (1943), 86--96.
 \bibitem{urabe87} T. Urabe. \emph{Elementary transformations of Dynkin graphs and singularities on quartic surfaces}. Invent. math. \textbf{87} (1987), 549--572.
\end{thebibliography}
\end{document}